\documentclass[twoside,11pt]{amsart}

\usepackage{amsmath,latexsym,amssymb,times,mathptm}

\def\demo{\noindent{\bf Proof. }}
\def\sqr#1#2{{\vcenter{\hrule height.#2pt
        \hbox{\vrule width.#2pt height#1pt \kern#1pt
                \vrule width.#2pt}
        \hrule height.#2pt}}}
\def\square{\mathchoice\sqr64\sqr64\sqr{4}3\sqr{3}3}
\def\QED{\hfill$\square$}

\def\p{{\mathfrak p}}
\def\q{{\mathfrak q}}
\def\a{{\mathfrak a}}
\def\m{{\mathfrak m}}

\def\frak{\mathfrak}
\def\J{\mathcal J}

\def\K{\mathcal K}

\def\MM{\mathcal M}
\def\A{\mathcal A}
\def\U{\mathcal U}

\def\cI{\check{I}}
\def\cH{\check{H}}

\newtheorem{Theorem}{Theorem}[section]
\newtheorem{Lemma}[Theorem]{Lemma}
\newtheorem{Corollary}[Theorem]{Corollary}
\newtheorem{Proposition}[Theorem]{Proposition}

\newtheorem{Notation and Discussion}[Theorem]{Notation and Discussion}
\newtheorem{Assumptions}[Theorem]{Assumptions}
\newtheorem{Remark}[Theorem]{Remark}
\newtheorem{Example}[Theorem]{Example}

\newcommand{\ov}{\overline}

\newcommand{\Z}{\mathbb{Z}}

\newcommand{\xvec}[1]{\ensuremath{x_{1},\ldots,x_{#1}}}
\newcommand{\lcm}{\ensuremath{\mathrm{lcm}}}

\newcommand{\core}{\ensuremath{\mathrm{core}}}
\newcommand{\adj}{\ensuremath{\mathrm{adj}}}
\newcommand{\np}{\ensuremath{\mathrm{NP}}}
\newcommand{\spec}{\ensuremath{\mathrm{Spec}}}

\newcommand{\mono}{\ensuremath{\mathrm{mono}}}
\newcommand{\Mono}{\ensuremath{\mathrm{Mono}}}
\newcommand{\aaa}{\ensuremath{\underline{\alpha}}}
\newcommand{\bbb}{\ensuremath{\underline{\beta}}}
\newcommand{\car}{\ensuremath{\mathrm{char}}}

\newcommand{\ms}{\medskip}

\begin{document}

\baselineskip=16pt

\title[The core of zero-dimensional monomial ideals]
{\Large\bf The core of zero-dimensional monomial ideals}

\author[C. Polini, B. Ulrich and M. Vitulli]
{ Claudia Polini, Bernd Ulrich \and Marie A. Vitulli}

\thanks{AMS 2000 {\em Mathematics Subject Classification}.
Primary 13B21; Secondary 13A30, 13B22, 13C40.}

\thanks{The first two authors were supported in part by the
NSF. The first author was also supported in part by the NSA}

\address{Department of Mathematics, University of Notre Dame,
Notre Dame, Indiana 46556} \email{cpolini@nd.edu}

\address{Department of Mathematics, Purdue University,
West Lafayette, Indiana 47907} \email{ulrich@math.purdue.edu}

\address{Department of Mathematics, University of Oregon,
Eugene, Oregon 97403} \email{vitulli@math.oregon.edu}

\keywords{Cores, monomial ideals, reductions,  Rees algebras,
canonical modules, adjoints, coefficient ideals.}

\vspace{-0.1in}

\begin{abstract}
The core of an ideal is the intersection of all its reductions. We
describe the core of a zero-dimensional monomial ideal $I$ as the
largest monomial ideal contained in a general reduction of $I$. This
provides a new interpretation of the core in the monomial case as
well as an efficient algorithm for computing it. We relate the core
to adjoints and first coefficient ideals, and in dimension two and
three we give explicit formulas.
\end{abstract}

\maketitle

\vspace{-0.2in}

\section{Introduction}
The purpose of this paper is to study the core of monomial ideals.
According to Northcott and Rees \cite{NR}, a subideal $J$ of an
ideal $I$ is a \textit{reduction} of $I$ provided $I^{r+1}= JI^r$
for some nonnegative integer $r$. In a Noetherian ring, $J$ is a
reduction of $I$ if and only if $I$ is integral over $J$.
Intuitively, a reduction of $I$ is a simplification of $I$ that
shares essential properties with the original ideal. Reductions are
highly non-unique, even minimal reductions (with respect to
inclusion) that are known to exist for ideals in Noetherian local
rings. Thus one considers the ${\it core}$ of the ideal $I$, written
$\core(I)$, which is the intersection of all reductions of $I$.

The core, introduced by Rees and Sally \cite{RS}, is in a sense
the opposite of the integral closure: the integral closure
$\overline{I}$ is the largest ideal integral over $I$, whereas
$\core(I)$ is the intersection of all ideals over which $I$ is
integral. The core appears naturally in the context of
Brian\c{c}on-Skoda theorems that compare the integral closure
filtration with the adic filtration of an ideal. It is also
connected to adjoints, multiplier ideals and coefficient ideals.

Huneke-Swanson, Corso-Polini-Ulrich,  Hyry-Smith, Polini-Ulrich, and
Huneke-Trung \cite{HS, CPU1, CPU2, HyS1, PU, HT, HyS2} gave explicit
formulas for cores in local rings (whose residue characteristic is
zero or large enough) by expressing them as colon ideals. For
certain classes of ideals, which include zero-dimensional ideals,
they showed that $\core(I) = J^{n+1}:I^n$, where $J$ is a minimal
reduction of $I$ and $n$ is sufficiently large. Moreover, Hyry and
Smith \cite{HyS1, HyS2} discovered an unforeseen relationship with
Kawamata's conjecture on the non-vanishing of sections of line
bundles. They proved that Kawamata's conjecture would follow from a
formula that essentially amounts to a graded analogue of the above
formula for the core.

The known formulas for the core usually require the ambient ring to
be \textit{local}. In contrast, in this paper we are primarily
interested in the core of 0-dimensional monomial ideals in
polynomial rings. Thus we start Section 2 by establishing the
expected colon formula for the core in the \text{global } setting,
for 0-dimensional ideals. For this we prove that the core of
0-dimensional ideals commutes with localization.

Let $R=k[x_1, \ldots, x_d]$ be a polynomial ring over an infinite
field $k$, write $\m=(x_1, \ldots, x_d)$, and let $I$ be a monomial
ideal, that is, an $R$-ideal generated by monomials. Even though
there may not exist any proper reduction of $I$ which is monomial
(or even homogeneous), the intersection of all reductions, the core,
is again a monomial ideal (because of the torus action, see for
instance \cite[5.1]{CPU1}). Lipman \cite{L} and Huneke-Swanson
\cite{HS} related the core to the adjoint ideal (see also \cite{Hy,
HyS1, HyS2, PU}). The integral closure and the adjoint of a monomial
ideal are again monomial ideals and can be described in terms of the
Newton polyhedron $\np(I)$ of $I$ \cite{How, HuS}. Such a
description cannot exist for the core, since the Newton polyhedron
only depends on the integral closure of the ideal, whereas the core
may change when passing from $I$ to $\overline{I}$. When attempting
to derive any kind of combinatorial description for the core of a
monomial ideal from the known colon formulas, one faces the problem
that the colon formulas involve non-monomial ideals, unless $I$ has
a reduction $J$ generated by a monomial regular sequence. Instead,
we exploit the existence of such non-monomial reductions to devise
an interpretation of the core in terms of monomial operations. This
is done in Section 3, where we prove that the core is the largest
monomial ideal contained in a `general locally minimal reduction' of
$I$.

Let $I$ be a $0$-dimensional monomial ideal in $k[x_1,\ldots,x_d]$
and $J$ an ideal generated by $d$ general $k$-linear combinations of
minimal monomial generators of $I$. Unless $I$ is generated by
monomials of the same degree, $J$ may not even be $\frak
m$--primary, but $J_{\m}$ is a minimal reduction of $I_{\m}$. Since
$I$ is $\m$-primary, there exist $n_i$ such that $x_i^{n_i}\in I$.
The regular sequence $\underline\alpha=x_1^{dn_1},\ldots,x_d^{dn_d}$
is contained in the core of $I_{\m}$ by the Brian\c con--Skoda
theorem. Hence $(J,\underline\alpha)_{\m}=J_{\m}$. Because
$K=(J,\underline\alpha)$ is a reduction of $I$ with $K_{\m} =
J_{\m}$, we call such $K$ a {\it general locally minimal reduction}
of $I$. As ${\rm core}(I)$ is a monomial ideal contained in $K$, it
is contained in ${\rm mono}(K)$, the largest monomial subideal of
$K$. In Theorem~\ref{core = mono} we actually show that $ \core(I) =
\mono(K)$. Notice that one cannot expect the inclusion $\core(I)
\subset \mono(K)$ to be an equality unless $K$ is far from being
monomial -- which is guaranteed by our general choice of $K$.

The idea behind the proof of Theorem~\ref{core = mono} is to show
that $\mono(K)$ is independent of the general locally minimal
reduction $K$. Using the inclusion reversing operation of linkage,
we express ${\rm mono}(K)$ in terms of
$\Mono((\underline\alpha)\colon K)$. Here $\Mono(L)$ denotes the
smallest monomial ideal containing an arbitrary ideal $L$, which can
be easily computed as it is generated by the monomial supports of
generators of $L$. We are able to show that
$\Mono((\underline\alpha)\colon K)$ does not depend on $K$, which
together with the equality ${\rm mono}(K) =
(\underline{\alpha}):\Mono((\underline\alpha)\colon K)$ gives the
independence of $\mono(K)$. The last equality is also interesting as
it establishes a link between $\mono$ and $\Mono$, and because it
yields an algorithm for computing $\mono$. A different algorithm can
be found in Saito-Sturmfels-Takayama \cite{SST}. Besides providing a
new, combinatorial interpretation of the core, the formula $
\core(I) = \mono(K)$ is in general more efficient computationally
than the colon formula ${\rm core}(I)=J^{n+1} \colon I^n$, as it
only requires taking colons of $d$-generated ideals. Furthermore the
new formula holds without any restriction on the characteristic.

Another way to find a combinatorial description of the core of a
monomial ideal is to express it as the adjoint of a power of the
ideal and use the known description of adjoints in terms of Newton
polyhedra. We pursue this approach in Section 4, where we show that
$\core(I)=\adj(I^d)$ if $I$ is a $0$-dimensional monomial ideal $I$
in a polynomial ring $k[x_1,\ldots,x_d]$ of characteristic zero and
all large powers of $I$ are integrally closed or nearly integrally
closed (see Theorem~\ref{adj2}, which uses Boutot's Theorem
\cite{B}, or Theorem~\ref{adj1} featuring a special case with an
elementary proof). On the other hand, the assumption on the integral
closedness is not always necessary, for in Sections 6 and 7 we
present classes of ideals in dimension two and three for which this
condition fails, whereas $\core( I) = \adj(  I^d)$. Our results of
Section 4 are based on the fact that both the core and the adjoint
can be related to components of the graded canonical module
$\omega_{R[It,t^{-1}]}$ of the extended Rees algebra $R[It,t^{-1}]$.
This approach also led us to study the core by means of the first
coefficient ideal $\cI$ of $I$. Let $D={\rm
End}(\omega_{R[It,t^{-1}]})$ denote the $S_2$-ification of the
extended Rees algebra of $I$ and define $\check{I}$ to be the
$R$-ideal with $D_1$ = $\cI t$; this ideal is also the first
coefficient ideal of $I$, the largest ideal that has the same zeroth
and first Hilbert coefficient as $I$  \, \cite{Shah, Ciu1}. As
remarked before, the core may change as one passes from $I$ to its
integral closure $\overline{I}$, however we show in
Theorem~\ref{checkI} that one can replace $I$ by any ideal between
$I$ and $\cI$ to compute the core, assuming that $I$ is a
0-dimensional monomial ideal in characteristic zero. If $I$ has a
reduction generated by a monomial regular sequence we prove in fact
that $\cI$ is the unique largest ideal integral over $I$ that shares
the same core (see Corollary~\ref{checkI3}).

In Sections 6 and 7 we determine explicitly the core of ideals
generated by monomials of the same degree, in a polynomial ring in
$d \leq 3$ variables. For instance, consider the case $d=2$ and
write $I=\mu (x^n, y^n, \{x^{n-k_i}y^{k_i}\})$ with $\mu$ a
monomial. We show that ${\rm core}(I)=\mu (x^{\delta},
y^{\delta})^{2 \frac{n}{\delta} -1}$ where $\delta={\rm
gcd}(\{k_i\}, n)$ (see Theorem~\ref{coredim2}). In particular if
$\mu=1$ and $\delta=1$, then the core of $I$ is a power of the
maximal ideal and ${\rm core}(I)$ equals ${\rm adj}(I^2)$ even
though $I$ need not be integrally closed (see
Corollary~\ref{coreadj2}).

\bigskip

\bigskip

\section{Preliminaries}

\medskip

In this section we prove some general facts about cores in rings
that are not necessarily local. First we deal with the behavior of
cores under localization. This issue was addressed in \cite{CPU1}
for local rings.  Now instead we assume that the ideal be
$0$-dimensional in order to assure that the core is a finite
intersection of reductions. We then use the results of \cite{PU,
HT, FPU} to obtain explicit formulas for the core in global rings.

\medskip

\begin{Proposition}\label{local global} Let $R$ be a Noetherian ring, $S$ a
multiplicative subset of $R$, and $I$ a $0$-dimensional ideal.
Then \[ \core(S^{-1}I) = S^{-1}\core(I).\]
\end{Proposition}

\begin{proof}
Notice that there exists an integer $N \ge 0$ such that $I^N \subset
J$ for every reduction $J$ of $I$ \, \cite[2.4]{Vas}. From this it
follows that  $\core(I)$ is 0-dimensional.  Hence $R/\core(I)$ is
Artinian, which implies that $\core(I)$ is a finite intersection of
reductions. Say $\core(I) = \bigcap_{i = 1}^t J_i$. The inclusion
$\core(S^{-1}I) \subset S^{-1}\core(I)$ follows from
\[
\core(S^{-1}I) \subset  \bigcap_{i = 1}^t  S^{-1}J_i =
S^{-1}\bigcap_{i = 1}^t  J_i = S^{-1}\core(I) .
\]

To prove that $S^{-1}\core(I) \subset \core(S^{-1}I)$  we will show
that every reduction of $S^{-1}I$ is the localization of a reduction
of $I$. Let $\J \subset S^{-1}R$ be a reduction of $S^{-1}I$ and
consider $J = \J \cap I$. Obviously $S^{-1}J = \J$. We claim that
$J$ is a reduction of $I$. It suffices to prove this locally at
every prime $\p$ of $R$. If $(\J \cap R)_{\p} =R_{\p}$ then
$J_{\p}=I_{\p}$. Now assume that $(\J \cap R)_{\p} \not= R_{\p}$.
For every minimal prime $\q$ of $\J \cap R$, the ideal $S^{-1}\q$ is
a minimal prime of $\J$, hence of $S^{-1}I$. Therefore $\q$ is a
minimal prime of $I$, showing that $\J \cap R$ is $0$-dimensional.
Hence $\p$ is a minimal prime of $\J \cap R$. Therefore as before
$S^{-1}\p$ is a minimal prime of $\J$, which gives
$R_{\p}=(S^{-1}R)_{S^{-1}\p}$. Hence $J_{\p}={\J}_{S^{-1}\p}$ is a
reduction of $I_{\p}$.
\end{proof}

\medskip

Let $R$ be a ring. Recall that if $J$ is a reduction of an $R$-ideal
$I$, then the {\it reduction number} $r_J(I)$ of $I$ with respect to
$J$ is the smallest nonnegative integer $r$ with $I^{r+1}=JI^r$. For
a sequence $\aaa=\alpha_1, \ldots, \alpha_s$ of elements in $R$ and
a positive integer $t$, we write $\aaa^t$ for the sequence
$\alpha_1^t, \ldots, \alpha_s^t$. If $L$ is a monomial ideal in a
polynomial ring with minimal monomial generators $\aaa=\alpha_1,
\ldots, \alpha_s$, write $L^{\langle t \rangle}=(\aaa^t)$.
\medskip

\begin{Lemma}\label{colon lemma}
Let $R$ be a Noetherian ring, and let $I$ be an
ideal with $g={\rm ht} \, I >0$ having a reduction generated by a
regular sequence $\aaa$. Then for  $t \ge r_{(\aaa)}(I)$ and $i\ge
0$,
\[(\aaa)^{t+i}: I^t= (\aaa^{t+i}) \colon I^{\, gt+(g-1)(i-1)} = (\aaa^{t+i}) \colon (I^{\, gt+(g-1)(i-1)},  \aaa^{t+i}).
\]
\end{Lemma}
\demo Since $\aaa$
 is  a regular sequence we have
\[
 (\aaa^{t+i}) \colon (\aaa)^{(g-1)(t+i-1)} = (\aaa)^{t+i}.
\]
Hence for $t \ge
r_{(\aaa)}(I)$,
\begin{eqnarray*}
(\aaa)^{t+i}:I^{\, t}& = & ((\aaa^{t+i}) \colon (\aaa)^{(g-1)(t+i-1)}):I^t  \\
&=&  (\aaa^{t+i}): (\aaa)^{(g-1)(t+i-1)} I^t  \\
&=& (\aaa^{t+i}) \colon I^{gt+(g-1)(i-1)} \\
&=& (\aaa^{t+i}) \colon (I^{gt+(g-1)(i-1)},  \aaa^{t+i}).
\end{eqnarray*}\QED

\medskip

We are now ready to state the formulas for the core that we will
use throughout:

\begin{Theorem}  \label{core formula 1} Let $R$ be a
Cohen-Macaulay  ring containing an infinite field $k$ and  $I$ a
0-dimensional ideal of height $d>0$
having a reduction generated by a regular sequence $\aaa$. Assume
that $\car \, k =0$ or $\car \, k > r_{(\aaa)}(I)$. Then for $t \ge
r_{(\aaa)}(I)$,
\[\core(I) =(\aaa)^{t+1} \colon I^t =   (\aaa^{t+1}) \colon I^{dt} = (\aaa^{t+1}) \colon (I^{dt}, \aaa^{t+1}).
\]
\end{Theorem}
\demo
 Proposition~\ref{local global}, \cite[3.7]{HT}, and
 \cite[3.4]{PU} show that
 $\core(I) = (\aaa)^{t+1}:I^t$ for $t \ge r_{(\aaa)}(I)$. The last
 two equalities follow from Lemma~\ref{colon lemma}.
\QED

\medskip

\begin{Remark}\label{kperfect}
{\rm If in Theorem~\ref{core formula 1} the ideal $I$ is unmixed
then the assumption
 that $I$ has a reduction generated by a regular
 sequence is automatically satisfied, as can be seen from basic element theory. For a more general result we refer to  \cite[Theorem]{Ly}.}
\end{Remark}

\medskip

In the graded case, the assumption on the characteristic in
Theorem~\ref{core formula 1} can be dropped:

\medskip

\begin{Theorem}  \label{core formula} Let $R$ be a Cohen-Macaulay
geometrically reduced positively graded ring over an infinite field
and I a $0$-dimensional
ideal of height $d>0$ generated by forms of the same degree.
Let $\aaa$ be a homogeneous regular sequence generating a
reduction of $I$. Then for $t \ge r_{(\aaa)}(I)$,
\[\core(I) =(\aaa)^{t+1} \colon I^t =   (\aaa^{t+1}) \colon I^{dt} = (\aaa^{t+1}) \colon (I^{dt}, \aaa^{t+1}).
\]
\end{Theorem}
\demo  By  \cite[4.1]{FPU} we have  $\core(I) = (\aaa)^{t+1}:I^t$
for $t \ge r_{(\aaa)}(I)$. The other two equalities follow from
Lemma~\ref{colon lemma}.\QED

\medskip

\begin{Remark}\label{kperfect2}
{\rm Notice that a regular sequence $\aaa$ as in Theorem~\ref{core
formula} always exists.}
\end{Remark}

\bigskip

\section{An algorithm}

\medskip

In this section we prove a formula for the core of $0$-dimensional
monomial ideals. This formula gives a new interpretation of the core
in terms of operations on monomial ideals and at the same time
provides an algorithm that is more efficient in general than the
formulas of Theorems~\ref{core formula 1} and \ref{core formula}.
Furthermore the new approach does not require any restriction on the
characteristic.

\medskip

\begin{Notation and Discussion}\label{mono}
{\rm Let  $R = k[\xvec{d}]$ be a polynomial ring over a field $k$.
For an $R$-ideal $L$ we let $\mono(L)$ denote the largest monomial
ideal contained in $L$ and $\Mono(L)$ the smallest monomial ideal
containing $L$. Note that $\Mono(L)$ is easy to compute, being the
ideal generated by the monomial supports of generators of $L$.  The
computation of $\mono(L)$ is also accessible; the algorithm provided
in \cite[4.4.2]{SST} computes $\mono(L)$ by multi-homogenizing $L$
with respect to $d$ new variables and then contracting back to the
ring $R$.  The ideal $\mono(L)$ can be computed in CoCoA with the
built-in command {\em MonsInIdeal} .}
\end{Notation and Discussion}

\medskip

From now on let $k$ be an infinite field and write $\m =
(\xvec{d})$ for the homogeneous maximal ideal of $R$. To begin we
will use linkage to give a new algorithm to compute $\mono(L)$ for
a class of ideals including $\m$-primary ideals.

\medskip

 \begin{Lemma} \label{linkage}  Let $L$ be an unmixed $R$-ideal
 of height $g$ and $\bbb \subset L$ a regular sequence consisting
 of $g$ monomials. Then
 \[\mono(L) = (\bbb): \Mono((\bbb):L).\]
\end{Lemma}
\demo
Notice that $(\bbb):\Mono((\bbb):L) \subset (\bbb):((\bbb):L)
\subset L$, where the last containment holds since $R/(\bbb)$ is
Gorenstein and $L$ is unmixed. Now observe that colons of monomial
ideals are monomial. Hence $(\bbb):\Mono((\bbb):L) \subset
\mono(L)$. The other inclusion follows from the following
containments.  First, $(\bbb):L \subset (\bbb):\mono(L)$. But
$(\bbb):\mono(L)$ is monomial and hence $\Mono((\bbb):L) \subset
(\bbb):\mono(L)$.
Therefore $\mono(L) \subset (\bbb):\Mono((\bbb):L)$. \QED

\medskip

\begin{Notation and Discussion}\label{minred}
{\rm  Now let $I$ denote an $\m$-primary monomial ideal.  For each
$i$ let $n_i$  be the smallest power of $x_i$ in $I$; such $n_i$
exist since $I$ is $\m$-primary. Write $\aaa = x_1^{dn_1}, \ldots ,
x_d^{dn_d}$ and let $J$ be an ideal generated by $d$ general
$k$-linear combinations of minimal monomial generators of $I$. If
the ideal $I$ is generated by forms of the same degree, $J$ is a
general minimal reduction of $I$ \, \cite[5.1]{NR}. In general
however, $I$ and $J$ may not even have the same radical.
Nevertheless, $J_{\m}$ is a general minimal reduction of $I_{\m}$ by
\cite[5.1]{NR}. Consider the ideal $K=(J,\underline\alpha)$. Observe
that the $\m$-primary ideal $K$ is a reduction of $I$. Thus ${\rm
core}(I) \subset {\rm mono}(K)$ since the core is a monomial ideal.
The Brian\c{c}on-Skoda theorem implies $(\aaa )_{\m} \subset
\core(I_{\m})$. Hence $K_{\m}=J_{\m}$, and whenever $I$ is generated
by forms of the same degree then $K=J$. We call $K$ a {\it general
locally minimal reduction} of $I$.}
\end{Notation and Discussion}

\medskip

In order to prove the equality $\core(I)= \mono(K)$ we need to
show that $\mono(K)$ is independent of $K$; by this we mean that
$\mono(K)$ is constant as the coefficient matrix defining $J$
varies in a suitable dense open set of an affine $k$-space:

\begin{Lemma}\label{independence} With assumptions as in \ref{mono} and in
\ref{minred}, the ideal $\Mono((\aaa):K)$ does not depend on the general
locally minimal reduction $K$.
\end{Lemma}
\demo
 Let $f_1, \ldots, f_n$ be minimal monomial generators of $I$.
 Let $\underline{z} = z_{i j}$, $1 \leq i \leq d$, $1 \leq j \leq n$,
 be variables and write
 $T = R[\underline{z}]$.
Let $\J$ denote the $T$-ideal generated by the $d$ generic linear
combinations $\sum_{j = 1}^n z_{i j}f_j$, $1 \leq i \leq d$, and let
$\K$ be the $T$-ideal  $(\J, \aaa)$. For $\underline{\lambda} =
\lambda_{ij}$, $1 \leq i \leq d$, $1 \leq j \leq n$, any elements in
$k$, we consider the maximal ideal $\MM = ({\mathfrak
m},\underline{z}-\underline{\lambda}) = ({\mathfrak m}, \{
z_{ij}-\lambda_{ij} \})$ of $T$. We identify the set $\A=\{\MM=
({\mathfrak m},\underline{z}-\underline{\lambda}) \,|\,
\underline{\lambda} \in k^{d n} \}$ with the set of $k$-rational
points of the affine space ${\mathbb A}_k^{d n}$. Write
$\pi_{\underline{\lambda}} : T \rightarrow R$ for the homomorphism
of $R$-algebras with $\pi_{\underline{\lambda}}(z_{ij}) =
\lambda_{ij}$. This map induces a local homomorphism $T_{\MM}
\rightarrow R_{\m}$, which we still denote by
$\pi_{\underline{\lambda}}$.

Notice that $\pi_{\underline{\lambda}}(\K) = K$ for
$\underline{\lambda}$ in a dense open subset $U_1 \subset {\mathbb
A}_k^{d n }$.


Now we claim that there is a dense open subset $U_2 \subset {\mathbb
A}_k^{d n }$ such that $\K_{\, \MM}$ is Cohen-Macaulay. Indeed, let
$N$ be a $(d-1)^{\rm st}$ syzygy of the $T$-ideal $\K$. The free
locus of $N$ is a dense open subset $\U$ of $\spec(T)$. It contains
$\m T$ since $N_{\m T}$ is a $(d-1)^{\rm st}$ syzygy of the ideal
$\K_{\m T}$ over the $d$-dimensional regular local ring $T_{\m T}$.
Thus intersecting $\U$ with $\A$ we obtain a dense open subset
 $U_2 \subset
{\mathbb A}_k^{d n }$ where $N_{\, \MM}$ is free. Since the ideal
$\K_{\, \MM}$ has height at least $d$ it is Cohen-Macaulay.

For every $\underline{\lambda} \in U_2$ the ideal $\K_{\MM}$ is
Cohen-Macaulay. Therefore $(\aaa):\K_{\MM}$ specializes according to
\cite[2.13]{HU}, that is,
$\pi_{\underline{\lambda}}((\aaa):\K_{\MM}) =
(\aaa):\pi_{\underline{\lambda}}(\K_{\MM}) $. Thus
$\pi_{\underline{\lambda}}((\aaa):\K)_{\m} =
((\aaa):\pi_{\underline{\lambda}}(\K))_{\m}$ because
$\pi_{\underline{\lambda}}(T_{\MM}) = R_{\m}$. On the other hand,
$\pi_{\underline{\lambda}} ((\aaa) \colon \K)$ is $\m$-primary since
$\aaa=\pi_{\underline{\lambda}}(\aaa) \subset
\pi_{\underline{\lambda}} ((\aaa) \colon \K)$.
Therefore $\pi_{\underline{\lambda}}((\aaa):\K) =
(\aaa):\pi_{\underline{\lambda}}(\K)$ for every $\underline{\lambda}
\in U_2$.


We think of $T$ as a polynomial ring in $\xvec{d}$ over
$k[\underline{z}]$.
Write the generators of $(\aaa):\K$ as sums of monomials in the
$x$'s with coefficients $g_1(z), \ldots , g_t(z)$. The $R$-ideal
$\Mono(\pi_{\underline{\lambda}}((\aaa) \colon \K))$ is
independent of $\underline{\lambda}$ for $\underline{\lambda} \in
U_3=D_{g_1 \cdots g_t}$.

For $\underline{\lambda} \in U_1 \cap U_2 \cap U_3$ the $R$-ideal
$K=\pi_{\underline{\lambda}}(\K)$  is a general locally minimal
reduction  of $I$ and $\Mono((\aaa) \colon K)=\Mono((\aaa) \colon
\pi_{\underline{\lambda}}(\K))=\Mono(\pi_{\underline{\lambda}}((\aaa)
\colon \K))$ does not depend on $\underline{\lambda}$.
\QED

\medskip

\begin{Corollary}\label{mono-ind} With assumptions as in \ref{mono} and in
\ref{minred}, the ideal $\mono(K)$ does not depend on the general
locally minimal reduction $K$. \end{Corollary} \demo The claim
follows from Lemmas~\ref{linkage} and \ref{independence}. \QED

\bigskip

We are now ready to prove the main result of this section.

\begin{Theorem}\label{core = mono} With assumptions as in \ref{mono} and in
\ref{minred},
 \[\core(I) = \mono(K) = (\aaa): \Mono((\aaa):K).
 \]
\end{Theorem}
\demo We already know that $\core(I) \subset \mono(K) $.
Furthermore $\mono(K) = (\aaa): \Mono((\aaa):K)$ by
Lemma~\ref{linkage}. Thus
 it suffices to show that $\mono(K) \subset \core(I)$.
From \cite[4.5]{CPU1}  it follows that
\begin{equation*} \core(I_{\m}) = (K_1)_{\m} \cap \ldots \cap (K_t)_{\m}
\end{equation*}
 for general locally minimal reductions $K_1, \ldots, K_t$ of $I$.
 According to Corollary~\ref{mono-ind} we may assume that $\mono(K)=\mono(K_i)$
 for $1\leq i \leq t$. Therefore $\mono(K) \subset K_1 \cap \ldots \cap K_t$ and thus
 $\mono(K)_{\m} \subset \core(I_{\m})= \core(I)_{\m}$, where the last
 equality holds by Proposition~\ref{local global}. Hence $\mono(K) \subset
 \core(I)$ as $\core(I)$ is $\m$-primary.
 \QED

\medskip

\begin{Remark}\label{example}
{\rm The above theorem gives a new interpretation of the core of a
monomial ideal $I$ as the largest monomial ideal contained in a
general locally minimal reduction of $I$. This idea can be easily
implemented in CoCoA using a script to obtain $d$ general elements
in the ideal $I$ and the built-in command {\em MonsInIdeal} to
compute $\mono(K)$.}
\end{Remark}

\begin{Remark}\label{remark-ex2}
{\rm The formula of Theorem~\ref{core formula 1} does not hold in
arbitrary characteristic (see \cite[4.9]{PU}). However, if $J$ and
$I$ are monomial ideals, $J^{n+1}:I^n$ is obviously independent of
the characteristic. On the other hand, the algorithm based on
Theorem~\ref{core = mono} works in any characteristic, but its
output, $\mono(K)$, is characteristic dependent. In fact we are now
going to exhibit a zero-dimensional monomial ideal $I$ for which
$\core(I)=\mono(K)$ varies with the characteristic. As $I$ has a
reduction $J$ generated by a monomial regular sequence this shows
that the formula of Theorem~\ref{core formula 1} fails to hold in
arbitrary characteristic even for $0$-dimensional monomial ideals.}
\end{Remark}

\medskip

\begin{Example} \label{Ex2}
{\rm  Let $R=k[x,y]$ be a polynomial ring over an infinite field
$k$, consider the ideal $I=(x^6, x^5y^3, x^4y^4, x^2y^8,  y^9)$, and
write $J=(x^6,y^9)$. One has $r_J(I)=2$. If $\car \,  k \not= 2$
then the formula of Theorem~\ref{core formula 1} as well as the
algorithm of Theorem~\ref{core = mono} give $\core(I)=J^3 \colon
I^2= J(x^4,  x^3y, x^2y^2,  xy^5, y^6)= (x^{10}, x^9y, x^8y^2,
x^7y^5, x^6y^6, x^4y^9,  x^3y^{10},  x^2y^{11}, xy^{14}, y^{15})$.
On the other hand, if $\car \, k =2$ then Theorem~\ref{core = mono}
shows that $\core(I)= (x^{10}, x^8y,  x^7y^5, x^6y^6, x^4y^9,
x^3y^{10}, x^2y^{11}, xy^{14}, y^{15}) \varsupsetneq J^3 \colon
I^2$. }
\end{Example}

\bigskip

\section{ The core, the first coefficient ideal and the adjoint}

\medskip

\begin{Notation and Discussion}\label{end}
{\rm Let $R$ be a Gorenstein ring, let $I$ be an $R$-ideal with
$g={\rm ht} \ I > 0$, and assume that $I$ has a reduction $J$ which
is locally a complete intersection of height $g$. Consider the
inclusions
\[A=R[Jt,t^{-1}] \subset B=R[It,
t^{-1}] \subset R[t,t^{-1}].\] Notice that $A$ is a Gorenstein ring.
We define $\omega_A = A t^{g-1} \subset R[t,t^{-1}]$ and write
${-}^{\vee}= {\rm Hom}_A(-, \omega_A)$, $F= {\rm Quot}(R[t])$. We
may choose $\omega_B= \omega_A :_{R[t,t^{-1}]} B =\omega_A :_F B
\simeq B^{\vee}$ as a graded canonical module of $B$. According to
 \cite[2.2.2]{PU},
 \begin{equation}\label{canonical}
 \omega_B = \oplus_i (J^{s+i-g+1} \colon I^s)t^i
\end{equation}
for every $s \ge r_J(I)$. Observe that $[\omega_B]_i=Rt^i$ for $i
\ll 0$. Write \begin{eqnarray*}
D &=& \omega_B :_{R[t,t^{-1}]} \omega_B \\
& = &  \omega_B :_F \omega_B \\
& = &  \omega_A :_F \omega_B \\
&=& A :_F (A :_F B) \\
&= & A :_{R[t,t^{-1}]} (A :_{R[t,t^{-1}]} B).
\end{eqnarray*}

\noindent Notice that $D \simeq {\rm End}_B(\omega_B) \simeq
B^{\vee\vee}$ is an $S_2$-fication of $B$. Define $\cI$ to be the
$R$-ideal with $[D]_1=\cI t$. One has $I \subset \cI \subset
\overline{I}$, and $\cI$ is the first coefficient ideal  of $I$ in
the sense of \cite{Shah, Ciu1, Ciu2}. Finally, write $C=R[\cI t,
t^{-1}]$. The inclusions $B\subset C \subset D$ are equalities
locally in codimension one in $A$, and hence upon applying $\omega_A
:_F - \simeq -^{\vee}$ yield equalities
\begin{equation}\label{eqS2}\omega_B=\omega_C=\omega_D.
\end{equation}  }
\end{Notation and Discussion}

\bigskip

\bigskip

We first give a formula expressing $D$ and $\cI$ in terms of colon
ideals. For this we need to consider an integer $u\ge 0$ such that
the graded canonical module of $B=R[It,t^{-1}]$ is generated in
degrees at most $g-1+u$ as a module over $A=R[Jt,t^{-1}]$. Whenever
$I$ is a monomial ideal one can take $u=0$, as we will see in
Theorem~\ref{comes out}. However, this is not longer true if $I$ is
not monomial and $B$ is not Cohen-Macaulay, see \cite[4.13]{PU}.

\begin{Theorem}\label{checkI2}
In addition to the assumptions of \ref{end} suppose that $R$ is
regular. Let $s \geq r_J(I)$ be an integer and $u \geq 0$ an integer
such that $J^{s+u+i} : I^s=J^i(J^{s+u}:I^s)$ for every $i \geq 0$.
One has
\[D = \oplus_i  (J^{i+u} \colon (J^{s+u} \colon I^s))t^i.\]
In particular
\[\check{I}= J^{1+u} \colon (J^{s+u}\colon I^s).\]
\end{Theorem}
\demo We need to prove that $
 D=A :_{R[t,t^{-1}]} (J^{s+u} :_R I^{s})t^u
$. The Brian\c{c}on-Skoda Theorem \cite[Theorem 1]{LS} gives
$I^{s+i} \subset J^{s+i-g+1}$ for every integer $i$, hence $J^i
\subset J^{s+i-g+1} :_R I^s$. Now  Equation~(\ref{canonical}) shows
that $A \subset \omega_B$. The same equation and our assumption also
give $[\omega_B]_i = (J^{s+u} :_R I^s)[t^u\omega_A]_i$ for $i \ge
g-1+u$. Hence writing $L=A+(J^{s+u} :_R I^s)t^u\omega_A$ we obtain
an exact sequence of graded $A$-modules
\[ 0 \longrightarrow L \longrightarrow \omega_B \longrightarrow N
\longrightarrow 0 \, ,
\]
with $N$ concentrated in finitely many degrees. It follows that $N$
has grade $\geq 2$.

Thus applying $\omega_A :_F - \simeq -^{\vee}$ yields
\begin{eqnarray*}
D & = & \omega_A :_F \omega_B \\
&=& \omega_A :_F L \\
&=& (\omega_A :_F A) \cap (\omega_A :_F (J^{s+u} :_R I^s)t^u\omega_A) \\
&=& \omega_A \cap (A :_F (J^{s+u} :_R I^s)t^u) \\
&=& A :_{\omega_A} (J^{s+u} :_R I^s)t^u.
\end{eqnarray*}
As $J^{g-1+u}\subset J^{s+u}:_R I^s$ we obtain
\[ J^{i+u} :_R (J^{s+u} :_R I^s) \subset J^{i+u} :_R
J^{g-1+u} = J^{i-g+1},
\]
where the last equality holds because ${\rm gr}_J(R)$ is
Cohen-Macaulay and ${\rm ht} \, J > 0$. Thus $A :_{R[t,t^{-1}]} (J^s
:_R I^s)t^u \subset \omega_A$, showing that
\[ A :_{R[t,t^{-1}]} (J^{s+u}
:_R I^s)t^u = A :_{\omega_A} (J^{s+u} :_R I^s)t^u= D. \]
 \QED

\medskip

In many cases all ideals between $I$ and $\cI$ have the same core:

\begin{Theorem}\label{checkI}
In addition to the assumptions of \ref{end} suppose that $R$
contains an infinite field $k$ with $\car \, k =0$ or $\car \, k >
r_J(I)$. Further assume that $R$ is local or $I$ is $0$-dimensional.
Then ${\rm core}(I)={\rm core}(\check{I})$.
\end{Theorem}
\demo By Proposition~\ref{local global} and \cite[4.8]{PU} we have
$J^{s+1} \colon \cI^s \subset \core(\cI)$ for $s \gg 0$. On the
other hand $\core(\cI) \subset \core(I)$ since $\cI$ is integral
over $I$. From Proposition~\ref{local global} and  \cite[4.5]{PU} we
obtain $\core(I)=J^{s+1} \colon I^s$. Finally,
Equations~(\ref{canonical}) and~(\ref{eqS2}) show that
\[(J^{s+1} \colon I^s)t^g =[\omega_B]_g=[\omega_C]_g =(J^{s+1} \colon
{\cI}^s)t^g.\]\QED

\medskip



\begin{Theorem} \label{check2} Let $R$ be a Gorenstein geometrically reduced positively
graded ring over an infinite field and I  a 0-dimensional
ideal generated by forms of the same degree.
Then ${\rm core}(I)={\rm core}(\check{I})$.
\end{Theorem}
\demo Let $J$ be a reduction of $I$ generated by a homogeneous
regular sequence and $s \gg 0$ an integer. As in the proof of
Theorem~\ref{checkI} one sees that $J^{s+1} \colon I^s =J^{s+1}
\colon {\cI}^s \subset \core(\cI) \subset \core(I)$. Furthermore
from Theorem~\ref{core formula} we obtain $\core(I)=J^{s+1} \colon
I^s$. \QED

\medskip
\begin{Assumptions}\label{monred}
{\rm Let $R=k[x_1, \ldots, x_d]$ be a polynomial ring over an
infinite field $k$ and write $\m = (\xvec{d})$ for the homogeneous
maximal ideal of $R$. Let $I \not= 0$ be  a  monomial ideal of
height $g$ and let $\a$ be an ideal generated by $g$ $k$-linear
combinations of the minimal monomial generators of $I$ . We assume
that $I$ has a reduction $J$ generated by a regular sequence of
monomials, and we write $r$ for the reduction number of $I$ with
respect to $J$. }
\end{Assumptions}
\medskip

Now our goal is to express $\cI$ as a colon ideal and to prove that
under certain conditions, $\cI$ is the unique largest ideal in
$\overline{I}$ having the same core as $I$. For this we need the
next theorem, which says that we may take $u=0$ in
Theorem~\ref{checkI2} provided we are in the setting of
\ref{monred}.

\begin{Theorem}\label{comes out} With assumptions as in \ref{monred}
one has for every $s \ge r$ and every $i \geq 0$,  \[J^{s+i} \colon
I^s=J^i(J^s \colon I^s)\] and
\[({\a}^{s+i} \colon I^s)_{\m}={\a}^i({\a}^s \colon
I^s)_{\m}.\]
\end{Theorem}
\demo  To prove the first equality write $f_1, \ldots, f_g$ for the
monomial generators of $J$. Clearly $J^i(J^s \colon I^s)\subset
J^{s+i} \colon I^s$. Notice also that $J^{s+i} \colon I^s \subset
J^{s+i} \colon J^s \subset J^i$ since $J$ is generated by a regular
sequence. Let $f$ be  a monomial contained in $J^{s+i} \colon I^s$,
and write $f=f_{j_{1}} \cdots f_{j_{i}} \cdot h$. Observe that
$f_{j_{1}} \cdots f_{j_{i}} \cdot h I^s=fI^s \subset J^{s+i}$.
Therefore $h I^s \subset J^{s+i} \colon (f_{j_{1}} \cdots
f_{j_{i}})= J^s$. Hence $h \in J^s \colon I^s$, which gives $f \in
J^i(J^s \colon I^s)$.

To prove the second equality notice that $r_{\a_{\m}}(I_{\m}) \leq
r$  \, \cite[3.4]{SUV}  and hence $({\a}^{s+i} \colon
I^s)_{\m}=(J^{s+i} \colon I^s)_{\m}$ by Equation~(\ref{canonical}).
Also observe that $(J^{s+i+1} \colon I^s)_{\m}=\a(J^{s+i} \colon
I^s)_{\m}$ whenever $i \ge i_0$ for some fixed integer $i_0$,
because $\omega_B \otimes_R R_{\m}$ is finitely generated as a
graded module over $R_{\m}[\a t, t^{-1}]$. Hence it suffices to
prove that $(J^{s+i} \colon I^s)_{\m}={\a}^i(J^s \colon I^s)_{\m}$
for each of the finitely many $i$ in the range $0 \le i \le i_0$.
We write $H =(J^{s+i} \colon I^s)_{\m}$ and $K=(J^{s} \colon
I^s)_{\m}$. Notice that $I^iK \subset H$ by
Equation~(\ref{canonical}) since $\omega_B$ is a $B$-module.

We complete $f_1, \ldots, f_g$ to monomial generators $f_1,
\ldots, f_n$ of $I$.
 Let $\underline{z} = z_{i j}$, $1 \leq i \leq g$, $1 \leq j \leq n$,
 be variables and write
 $T = R_{\m}[\underline{z}]$.
Let $\J$ denote the $T$-ideal generated by the $g$ generic linear
combinations $\sum_{j = 1}^n z_{i j}f_j$, $1 \leq i \leq g$.
Notice that $\J^i KT \subset HT$ as $\J \subset IT$. Since $H =J^i
K$ and $\J$ specializes to $J_{\m}$ modulo $( \{
z_{ij}-\delta_{ij} \})$, it follows that $HT = \J^iKT + [( \{
z_{ij}-\delta_{ij} \}) \cap HT] $. Consider the maximal ideal $\MM
= ({\mathfrak m},\underline{z}-\underline{\delta}) = ({\mathfrak
m}, \{ z_{ij}-\delta_{ij} \})$ of $T$. As
$\underline{z}-\underline{\delta}$  form a regular sequence on
$T_{\MM}$ and $T_{\MM}/HT_{\MM}$, we conclude that $HT_{\MM} =
\J^iKT_{\MM}$ according to Nakayama's Lemma.
For $\underline{\lambda} = \lambda_{ij}$, $1 \leq i \leq g$, $1
\leq j \leq n$, any elements in $k$, we consider the maximal ideal
$\MM_{\lambda} = ({\mathfrak m},\underline{z}-\underline{\lambda})
= ({\mathfrak m}, \{ z_{ij}-\lambda_{ij} \})$ of $T$. We identify
the set $\A=\{\MM_{\lambda} \,|\, \underline{\lambda} \in k^{g n}
\}$ with the set of $k$-rational points of the affine space
${\mathbb A}_k^{g n}$.
Since the two ideals $HT$ and $\J^iKT$ coincide locally at
$\MM=\MM_{\delta}$ the same holds locally at $\MM_{\lambda}$ for
$\underline{\lambda}$ in a dense open neighborhood of
$\underline{\delta}$ in $ {\mathbb A}_k^{g n }$. Specializing modulo
$\underline{z}-\underline{\lambda}$ we conclude that $H= \a^i K
$.\QED

\medskip

\begin{Corollary}\label{check-as-colon}
 With assumptions as in \ref{monred} one has for every $s \geq r$,
\[\check{I}= J \colon (J^s\colon I^s)\]
and
 \[\cI_{\m}= \a_{\m} \colon
(\a_{\m}^s\colon I_{\m}^s).\]
\end{Corollary}
\demo We use Theorems~\ref{checkI2} and \ref{comes out}. \QED
\medskip

\begin{Corollary}\label{checkI4} In addition to the assumptions of
\ref{monred} let $H$ be an ideal integral over $I$. If $J^{t+i}
\colon H^t= J^{t+i} \colon I^t$ for some $i \geq 0$ and $t \gg 0$,
then $\omega_{R[Ht,t^{-1}]}=\omega_{R[It,t^{-1}]}$.
\end{Corollary}
\demo Write $A=R[Jt,t^{-1}]$. We have an inclusion of finitely
generated graded $A$-modules
\[ \omega_{R[Ht,t^{-1}]} \subset \omega_{R[It,t^{-1}]}.
\]
By our assumption these modules coincide in degree $g+ i -1$
according to Equation~(\ref{canonical}). By Theorem~\ref{comes out}
the canonical module $\omega_{R[It,t^{-1}]}$ is generated in degrees
$ \leq g-1$ as an $A$-module, which forces the two modules to be the
same in degrees $ \geq g + i -1$. Furthermore the two modules
coincide in degrees $ \ll 0$. Since they satisfy $S_2$ it then
follows that they are equal. \QED

\ms

\begin{Corollary}\label{checkI3} In addition to the assumptions of
\ref{monred} suppose that $I$ is $0$-dimensional.
\begin{itemize}
\item[(a)] Let $H$ be an ideal integral over $I$ with the same
core as $I$. If $H$ and $I$ are generated by forms of the same
degree or if ${\rm char}\, k=0$, then
$\omega_{R[Ht,t^{-1}]}=\omega_{R[It,t^{-1}]}$.

\item[(b)] If ${\rm char}\, k=0$ then
the ideal $\check{I}$ is the largest ideal integral over $I$ with
the same core as $I$.
\end{itemize}
\end{Corollary}
\demo To prove part (a) notice that $J^{t+1} \colon
I^t=\core(I)=\core(H)=J^{t+1} \colon H^t$ for $t \gg 0$ by the first
equality in Theorems~\ref{core formula} or \ref{core formula 1}. Now
apply Corollary~\ref{checkI4}.

Part (b) follows from part (a). Indeed, by (a) if
$H$ is an ideal integral over $I$ with the same core as $I$ then
$\cI=\cH$. On the other hand, $\core(I)= \core(\cI)$ by
Theorem~\ref{checkI}. \QED

\medskip

The next corollary shows that in some cases the Rees ring of a
monomial ideal is Cohen-Macaulay if it satisfies $S_2$. Monomial
algebras in general are Cohen-Macaulay provided they are normal, but
the $S_2$ property does not suffice  \cite[Theorem 1 and Remark
4]{H}.

\begin{Corollary} In addition to the assumptions of
\ref{monred} suppose that $d=2$. One has:
\begin{itemize}
\item[(a)] $r_{J}(\check{I}) \leq 1$.

\item[(b)] $R[\check{I}t]$ is the $S_2-{\rm ification \ of \ }
R[It]$ and it is Cohen-Macaulay.

\item[(c)]If $R[It]$ satisfies $S_2$ then it is Cohen-Macaulay.
\end{itemize}
\end{Corollary}
\demo To prove part (a) we may replace $I$ by $\cI$ to assume $\cI =
I$. Observe that by Corollary~\ref{check-as-colon}, $I_{\m}= \a_{\m}
\colon (\a_{\m}^s\colon I_{\m}^s)$ for $s \gg 0$. However, $\a_{\m}
\subset \a_{\m}^s \colon I_{\m}^s$ according to the
Brian\c{c}on-Skoda Theorem \cite[Theorem 1]{LS}.   Therefore
$\a_{\m} \colon I_{\m} = \a_{\m} \colon (\a_{\m} \colon (\a_{\m}^s
\colon I_{\m}^s))= \a_{\m}^s \colon I_{\m}^s$. Since $\a_{\m}^s
\colon I_{\m}^s$ is the degree $g-1$ component of the canonical
module of $R_{\m}[I_{\m}t,t^{-1}]$, it does not depend on $\a_{\m}$.
Hence the ideal $I_{\m}$ is balanced \cite[3.6]{U2}. Therefore,
$I_{\m}$ has reduction number at most $1$ according to
\cite[4.8]{U2}. It follows that $r_{J}(I) \leq 1$.

To prove (b) and (c) observe that part (a), \cite[3.1]{VV}, and
\cite[3.10]{GS} imply the Cohen-Macaulayness of the Rees algebra of
$\cI_{\m}$ and hence of $\cI$. \QED

\medskip

We now turn to the relationship between cores and adjoints as
defined in \cite[1.1]{L}. Whenever the core is an adjoint one has a
combinatorial description of the former in terms of a Newton
polyhedron. In fact Howald has shown that if $I$ is a monomial ideal
then its adjoint (or multiplier ideal) $\adj(I)$ is the monomial
ideal with exponent set $\{ \alpha \in \Z_{\ge 0}^d \mid \alpha +
\mathbf{1} \in
 \np^{\circ}(I) \}$, where $\mathbf{1} = (1,1, \ldots, 1) \in \Z_{\ge 0}^d$ and
 $\np^{\circ}(I)$ denotes the interior of the Newton polyhedron of $I$ \, \cite[Main Theorem]{How} (see also \cite[16.5.3]{HS2}).

\begin{Theorem} \label{adj1}  In addition to the assumptions of
\ref{monred} suppose that $I$ is $0$-dimensional.
 Assume that ${\rm char} \, k
=0$, $\car \, k > r_J(I)$, or $I$ is generated by monomials of the
same degree.
 If $\overline{I^{dt}} \subset (I^{dt}, J^{\langle t+1
\rangle})$ for some  $t \ge \max\{ r_J(I), d-1\}$, then $\core(I)
= \mathrm{adj}(I^d)$.
\end{Theorem}

\demo One has $\mathrm{adj}(I^d) \subset \mathrm{adj}(I_{\m}^d) \cap
R$ by the definition of the adjoint. On the other hand
\cite[1.4.1(ii)]{L} shows that  $\mathrm{adj}(I_{\m}^d) \subset
\core(I_{\m})$. Finally $\core(I_{\m}) \cap R = \core(I)$ according
to Proposition~\ref{local global}. Therefore $\mathrm{adj}(I^d)
\subset \core(I).$

To show the reverse inclusion notice that  $\core(I) = J^{\langle
t+1 \rangle} : I^{dt} =  J^{\langle t+1 \rangle} :
\overline{I^{dt}}$, where the first equality holds by
Theorems~\ref{core formula 1} and \ref{core formula}, and the second
equality follows from our assumption on $I$. Thus it suffices to
show that $ J^{\langle t+1 \rangle} : \overline{I^{dt}} \subset
\mathrm{adj}(I^d).$

Write $J = (x_1^{n_1}, \ldots , x_d^{n_d})$ and $L = \lcm(n_1,
\ldots, n_d)$. Consider the vectors $\mathbf{n}=(n_1, \ldots,
n_d)$, $\mathbf{\omega} = (L/n_1, \ldots, L/n_d)$ and
$\mathbf{1}=(1, \ldots, 1)$ in $\Z_{\ge 0}^d.$
Let  $x^{\alpha} \not\in \mathrm{adj}(I^d).$ We need to show that
$x^{\alpha} \not\in J^{\langle t+1 \rangle} : \overline{I^{dt}}.$ As
$J^{\langle t+1 \rangle} \subset J^d \subset \mathrm{adj}(I^d)$ we
conclude $x^{\alpha} \not\in J^{\langle t+1 \rangle}.$ Thus writing
$\beta= (t+1)\mathbf{n} - \alpha - \mathbf{1}$, we have $\beta \in
\Z_{\ge 0}^d$ and  $x^{\alpha}x^{\beta} \not\in J^{\langle t+1
\rangle}$. It remains to prove that $x^{\beta} \in
\overline{I^{dt}}=\overline{J^{dt}}$ or equivalently that $\omega
\cdot \beta \ge dtL.$ Indeed,  as $x^{\alpha} \not\in
\mathrm{adj}(I^d)=\mathrm{adj}(J^d),$ \cite[Main Theorem]{How} (see
also \cite[16.5.3]{HS2}) gives $\omega \cdot \alpha \le dL - \omega
\cdot \mathbf{1}$.  Hence
\begin{eqnarray*}
\omega \cdot \beta &=&  (t+1) \omega \cdot
\mathbf{n} - \omega \cdot \alpha - \omega \cdot \mathbf{1}\\
&=& (t+1)dL - \omega \cdot \alpha - \omega \cdot \mathbf{1} \\
& \ge & (t+1)dL - (dL - \omega \cdot \mathbf{1})- \omega \cdot
\mathbf{1}\\ &=& dtL \,.
\end{eqnarray*} \QED

\medskip
In characteristic $0$ one has a characterization for when
$\core(I)=\adj(I^d)$ even when the monomial ideal $I$ does not
have a reduction generated by a regular sequence of monomial.
However, the proof of this fact, which generalizes
\cite[5.3.4]{HyS1}, is less elementary than the one above.

\medskip

\begin{Theorem} \label{adj2} Let $R=k[x_1, \ldots, x_d]$ be a polynomial ring
over a field $k$ of characteristic $0$. Let $I$ be  a
$0$-dimensional monomial ideal and let $\aaa$ be a regular sequence
generating a reduction of $I$. Then \[\adj(I^d)= (\aaa)^{t+1} :
\overline{I^{t}} \subset (\aaa)^{t+1} : I^{t} =\core(I)\] for every
$t \geq {\rm max} \{ r_{(\aaa)}(I), \,  d-1\}$, and equality holds
if and only if $\overline{I^{dt}} \subset (I^{dt}, \aaa^{t+1})$ for
some  $t \ge \max\{ r_{(\aaa)}(I), d-1\}$.
\end{Theorem}
\demo Let $\overline{B}$ denote the integral closure  of $B=R[It,
t^{-1}]$ in $R[t,t^{-1}]$. According to \cite[Proposition 1]{H} the
integral closure ${\ov B}$ is a direct summand of a polynomial ring
over $k$, hence \cite[Th$\acute{{\rm e}}$or$\grave{{\rm e}}$me]{B}
shows that $\overline{B}$ has only rational singularities. Likewise
$\overline{R[It]}$ is Cohen-Macaulay by the same references or
\cite[Theorem 1]{H}. According to Proposition~\ref{local global} and
since $\adj(I^d)= \cap \, \adj(I_{\m}^d)$, where the intersection is
taken over all maximal ideals $\m$ of $R$, we may replace $R$ by any
of its localizations $R_{\m}$. As $\overline{B}$ has rational
singularities, one obtains $\adj(I^d)= [\omega_{\overline{B}}]_{d}$,
which can be deduced from \cite[1.3.1]{L} (see \cite{U} for
details). According to \cite{PUVV} the Cohen-Macaulayness of
$\overline{R[It]}$ implies that
$\overline{I^j}=(\aaa)^{j-d+1}\overline{I^{d-1}}$ for every $j\geq
d-1$. Now a computation as in \cite[2.2.2]{PU} yields
$[\omega_{\overline{B}}]_d=  (\aaa)^{t+1} : \overline{I^{t}}=
(\aaa^{ t+1}) : \overline{I^{dt}}$ for every $t \ge d-1$, where the
last equality follows as in Lemma~\ref{colon lemma}. Therefore
$\adj(I^d)=(\aaa)^{t+1} : \overline{I^{t}} =(\aaa^{ t+1}) :
\overline{I^{dt}}$. On the other hand $\core(I)=(\aaa)^{t+1} : I^{t}
=(\aaa^{t+1}) : I^{dt} $ for every $t \ge  r_{(\aaa)}(I)$ according
to Theorem~\ref{core formula 1}, and the assertion follows. \QED

\medskip

Notice that if equality holds in the previous theorem then
$\core(I)=\core(\overline{I})$. This condition is necessary for the
core to be the adjoint of $I^d$ as $\adj(I^d)
=\adj(\overline{I}^d)\subset \core(\overline{I})\subset \core(I)$.
On the other hand, the next example shows that the core may not
coincide with the adjoint even if the monomial ideal $I$ is
integrally closed.


\begin{Example} \label{Ex3} {\rm  Let  $k[x,y,z]$ be a polynomial
ring over an infinite field $k$ with $\car \, k \not= 2$ and let
$\m$ denote the homogeneous maximal ideal. Consider the ideal
$I=\overline{(x^3,y^4,z^5)}$ and write $J=(x^3,y^4,z^5)$. One has
$r_J(I)=2$. From the formula of Theorem~\ref{core formula 1} we
obtain $\core(I)=\m I^2$. Notice that $x^2y^3z^4 \not\in \m I^2$,
whereas $(x^2y^3z^4)^2 \in (\m I^2)^2$. Thus $\core(I)$ is not
integrally closed although $I$ is. In particular $\core(I)$ cannot
be an adjoint ideal because adjoints are always integrally closed.
Also notice that the Rees algebra $R[It]$ is Cohen-Macaulay because
$I$ is integrally closed with $r_J(I)\leq 2$, see \cite[p. 317]{Hu},
\cite[Theorem 1]{I}, \cite[3.1]{VV}, \cite[3.10]{GS}}.
\end{Example}

\bigskip

\section{ The core in weighted polynomial rings}

\medskip
For a positively graded ring $S$ and a positive integer $n$ we let
$S_{\geq n}$ denote the homogeneous $S$-ideal $\oplus_{\ell \geq n}
S_{\ell}$. Notice that $S_{\geq n}$ is not necessarily generated in
degree $n$. In this section we study the core of ideals of the form
$S_{\geq n}$, where $S$ is a weighted polynomial ring. The case of
section rings of line bundles has been been considered by Hyry and
Smith in connection with a conjecture by Kawamata (see \cite{HyS1,
HyS2}). For us, the ideals $S_{\geq n}$ are mainly interesting
because they shed light on the core of monomial ideals in standard
graded polynomial rings, as will be explained in
Section~\ref{secdim3}.

\medskip

\begin{Lemma}\label{Sn}
Let $R=k[x_1, \ldots, x_d]$ be a polynomial ring over a field $k$,
$S=k[x_1^{a_1}, \ldots, x_d^{a_d}]$,  $n$ a multiple of $\, {\rm
lcm}(a_1, \ldots, a_d)$, and $J$ the $S$-ideal generated by
$x_1^{n}, \ldots, x_d^{n}$. The following hold:
\begin{itemize}
\item[(a)] $J^i$ is a reduction of $S_{\geq in}$ for every $i \geq
1$.

\item[(b)]If the $S$-ideal $S_{\geq n}$ is normal then
\[J^{\langle t+1 \rangle} :_{S} (S_{\ge n})^{dt}=J^{\langle t+1
\rangle} :_{S} S_{\geq dnt}=S_{\geq dn-\sum a_i +1} \ \ \ {\rm for \
} t \ge d-1.\]
\end{itemize}
\end{Lemma}
\demo  For every monomial  $f \in S_{\geq in}$ we have $f^n \in
J^{in}$. This gives part (a).

To prove part (b) notice that  $(S_{\geq n})^{dt}= S_{\geq dnt}$
 by part (a) as $(S_{\geq n})^{dt}$ is integrally closed.  Thus it
suffices to show the second equality. Since $t \ge d-1$ we have
$J^{\langle t+1 \rangle} \subset S_{\geq n(t+1)} \subset S_{\geq
dn-\sum a_i +1}$, and we may pass to the ring $A=S/J^{\langle t+1
\rangle}$. Notice that $A$ is an Artinian graded Gorenstein ring
with socle degree $dn(t+1)-\sum a_i$.
Therefore $0 \colon_A (A_{\geq dnt})=A_{\geq dn-\sum a_i +1}$.
Indeed, to see that the left hand side is contained in the right
hand side, let $f \not= 0$ be a homogeneous element in $ 0
\colon_A (A_{\geq dnt})$. There exists a homogeneous element
$\lambda \in A$ such that $ 0 \not= \lambda f \in {\rm soc}(A)$.
In particular ${\rm deg}(\lambda) < dnt$ and ${\rm deg}(\lambda f)
= dn(t+1)-\sum a_i $.  This implies ${\deg}(f) \geq dn-\sum a_i
+1$, hence $f \in A_{\geq dn-\sum a_i +1}.$\QED

\medskip

\begin{Proposition}\label{PcoreSn} Let $R=k[x_1, \ldots, x_d]$ be a polynomial ring
 over an infinite field $k$, $S=k[x_1^{a_1},
\ldots, x_d^{a_d}]$, and  $n$ a multiple of  ${\rm lcm}(a_1, \ldots,
a_d)$. Assume that $\mathrm{char } \, k =0$ or the $S$-ideal
$S_{\geq n}$ is generated by monomials of degree $n$. If $\,S_{\geq
n}$ is a normal $S$-ideal then ${\rm core}(S_{\geq n}) =S_{\geq
dn-\sum a_i +1}$.
\end{Proposition}
\demo The assertion follows from Theorems~\ref{core formula 1} and
\ref{core formula}, and Lemma~\ref{Sn}. \QED



\medskip

\begin{Corollary}\label{CcoreSn} Let $R=k[x_1, \ldots, x_d]$ be a polynomial ring
 over an infinite field $k$, $S=k[x_1^{a_1},
\ldots, x_d^{a_d}]$, $a={\rm lcm}(a_1, \ldots, a_d)$, and $n=s a$.
Assume that ${\rm char} \, k =0$ or the $S$-ideal $S_{\geq n}$ is
generated by monomials of degree $n$. If $s \geq d-1$ then ${\rm
core}(S_{\geq n}) =S_{\geq dn-\sum a_i +1}$.
\end{Corollary}
\demo By \cite[3.5]{RRV} the $S$-ideal $S_{\geq n}$ is normal. Now
the assertion follows from Proposition~\ref{PcoreSn}. \QED

\medskip

\begin{Corollary}\label{coreSn}
Let $k[x,y,z]$ be a polynomial ring  over an infinite field $k$,
$S=k[x^a, y^b, z^c]$ with $a,b,c$ pairwise relatively prime, and
$n$ a multiple of $abc$. Assume that ${\rm char} \, k =0$ or the
$S$-ideal $S_{\geq n}$ is generated by monomials of degree $n$.
Then ${\rm core}(S_{\geq n}) =S_{\geq 3n-a-b-c +1}$.
\end{Corollary} \demo
The $S$-ideal $S_{\geq n}$ is normal according to \cite[3.13]{Vi}
and \cite[3.5]{RRV}. Again the assertion follows from
Proposition~\ref{PcoreSn}. \QED

\medskip
The next example shows that Proposition~\ref{PcoreSn} does not
hold without the normality assumption.

\begin{Example}{\rm Let $k[x,y,z]$ be a polynomial ring over a
field $k$ with $\car \, k =0$ and consider the subring
$S=k[x^{30},y^{35},z^{42}]$. We take  $n={\rm lcm}(30,35,42)=210$,
in which case $3n-a-b-c +1=524$. It turns out that $S_{\geq 524}
\subsetneq \core(S_{\geq 210}) \subsetneq \overline{\core(S_{\geq
210})}= S_{\geq 520}$. }
\end{Example}

\bigskip

\section{Monomials of the Same Degree: Dimension two}

\medskip

In this section we prove a formula for the core of ideals
generated by monomials of the same degree in a polynomial ring in
two variables. We start with a number theoretic lemma.

\begin{Lemma}\label{numeric}  Let $ k_1 , \ldots ,  k_s $ be non negative integers,
$n$ a positive integer, and write $\delta =  \gcd(k_1, \ldots, k_s,
n)$. Every integer $t$ divisible by $\delta$ can be written in the
form
\begin{equation*}
t = \alpha n + \sum_{i = 1}^s \beta_ik_i \, ,
\end{equation*}
where  $\beta_i \ge 0$ for all $i$ and $\sum_{i = 1}^s \beta_i <
n/\delta.$  Furthermore, if $t \gg 0$ we can take $\alpha \ge 0$.

\end{Lemma}
\demo The second assertion follows trivially from the first, since
$\sum \beta_i < n/\delta$ and $n$ and the $k_i$ are fixed.

Replacing $t$, $k_i$, $n$ by $t/\delta$, $k_i/\delta$, and
$n/\delta$, respectively, we may assume that $\delta = 1$. For any
$t \in \Z$, we can write $t = \alpha n + \sum_{i = 1}^s \beta_ik_i$
where $\alpha, \beta_i \in \Z$ since $\gcd(k_1, \ldots, k_s, n) =1$.
We proceed by induction on $s$.  Let $s = 1$. Write $\beta_1 = qn +
r$ with $0 \le r \le n-1$.  Then $t = \alpha n + \beta_1 k_1 =
(\alpha + qk_1)n + rk_1$.   So the  assertion holds for $s = 1$.

Now assume $s > 1$ and the first assertion holds for $s - 1$. Let
$\delta_j = \gcd(k_1, \ldots , \widehat{k_j} , \ldots , k_s, n)$ for
$1 \leq j \leq s$.  If $\delta_j = 1$ for some $j$ then the
conclusion follows from the induction hypothesis.  So assume that
$\delta_j > 1$ for all $j$. For each $1 \leq j \leq s$ choose a
prime $p_j$ that divides $\delta_j$; notice that $p_j \nmid k_j$.
Hence $p_1, \ldots , p_s$ are distinct primes, $\prod p_j \mid n$
and $\prod_{j \ne i} p_j \mid k_i$.  Thus $ \prod_{j \ne i} p_j \mid
\gcd(n,k_i)$ and $ \prod_{j \ne i} p_j \ge 2^{s-1} \ge s$, hence $
\gcd(n,k_i) \ge s$.
Changing $\beta_i$ modulo $n/\gcd(n,k_i)$ using the division
algorithm, we can assume that $0 \le \beta_i \le
\frac{n}{\gcd(n,k_i)} - 1 \le \frac{n}{s} - 1$ and hence $\sum
\beta_i \le n-1$. \QED

\medskip

\begin{Assumptions}\label{dim2}
{\rm Let $R=k[x,y]$ be a polynomial ring over a field $k$ and write
$\m$ for the homogeneous maximal ideal of $R$. Let $I$ be an
$R$-ideal generated by monomials of the same degree. Write $I = \mu
(x^n, y^n, x^{n-k_1}y^{k_1}, \ldots, x^{n-k_s}y^{k_s})$ with $\mu$ a
monomial and $0 < k_1 < \cdots < k_s < n$, and set $\delta =
\gcd(k_1, \ldots, k_s, n)$. }
\end{Assumptions}

\ms

\begin{Lemma}\label{t large}  In addition to the
assumptions of \ref{dim2} suppose that $\mu=1$ and $\delta=1$. Then
for $t \gg 0$,
\[
\m^{2nt} \subset I^{2t} + (x^{n(t+1)}, y^{n(t+1)}). \]
\end{Lemma}
\demo Consider a monomial generator $x^uy^v$ of $ \m^{2nt}$. Thus
$u+v = 2nt$ and we may assume $u < n(t+1)$ and $v < n(t+1)$.  Since
$u + v = 2nt = n(t+1) + n(t-1)$, we must have $v > n(t-1)$.  By
Lemma~\ref{numeric} we can write
$$v = \alpha n + \sum_{i = 1}^s \beta_ik_i,$$
where $ \beta_i \ge 0$ and $\sum_{i = 1}^s \beta_i \le n-1$. As $v
> n(t-1)$ and $t \gg 0$, we can take $\alpha \ge  0$; we also have
$\alpha \le  t$ since $v < n(t+1)$.

Now
\begin{eqnarray*}
u &= & 2nt - \alpha n - \sum \beta_ik_i \\
 & = & 2nt - \alpha n - \sum \beta_in + \sum \beta_i(n - k_i) \\
  & = & (2 t - \alpha   - \sum \beta_i)n + \sum \beta_i(n - k_i).
\end{eqnarray*}
Notice that $2 t - \alpha   - \sum \beta_i \ge 0$, because $t \gg 0$
and $\alpha + \sum \beta_i \le t + n-1 \le 2t$.  Thus
\[ (u,v) =  (2 t - \alpha   - \sum \beta_i)(n,0) +  \sum \beta_i(n - k_i,k_i) + \alpha(0,n) \]
is the exponent of a monomial in $I^{2t}$. \QED

\medskip

We are now ready to prove the main theorem of the section.

\begin{Theorem}\label{coredim2}   In addition to the
assumptions of \ref{dim2} suppose that $k$ is an infinite field.
Then
$$\core(I) = \mu \, (x^\delta,y^\delta)^{2\frac{n}{\delta} - 1}.$$
\end{Theorem}
\demo    First, we may assume $\mu=1$, since ${\rm
  core}(\mu \, I)= \mu \ {\rm core}(I)$ for any non zero divisor $\mu$.
  Passing to the subring $k[x^\delta,y^\delta]$ over which $k[x,y]$ is flat,
we may further suppose that $\delta=1$.   Indeed the core of
$0$-dimensional ideals is preserved by flat base change according to
Proposition~\ref{local global} and \cite[4.8]{CPU1}. Now we are left
to prove that ${\rm core}(I)= {\mathfrak m}^{2n-1}$. But
\begin{eqnarray*}
\core(I) &=& (x^{n(t+1)}, y^{n(t+1)}):(I^{2t},x^{n(t+1)}, y^{n(t+1)})
\hspace{.7cm}   \mbox{by Theorem~\ref{core formula}} \\
& = & (x^{n(t+1)}, y^{n(t+1)}):\m^{2nt}    \hspace{3cm}
\mbox{by Lemma~\ref{t large} }\\
& = &  \m^{2n-1}.
\end{eqnarray*}\QED

\ms

\begin{Corollary}
In addition to the assumptions of \ref{dim2} suppose that $\mu=1$
and $\delta=1$. Then $\check{I} =\m^n$.
\end{Corollary}
\demo We may assume that $k$ is infinite. By Theorem~\ref{coredim2}
we have $\core(I)= \core(\m^n)$. Now the assertion follows from
Corollary~\ref{checkI3}(a). \QED

\medskip

For any integrally closed ideal $I$ in a two-dimensional regular
local ring it is known that $\core(I)=\adj(I^2)$, by work of Huneke
and Swanson and of Lipman \cite{HS, L}. The next corollary shows
that this equality may hold even for ideals that are far from being
integrally closed.

\begin{Corollary}\label{coreadj2}
In addition to the assumptions of \ref{dim2} suppose that $k$ is an
infinite field, $\mu=1$, and $\delta=1$. Then $\core(I)={\rm
adj}(I^2)$.
\end{Corollary}
\demo The assertion follows from Theorem~\ref{adj1} via Lemma~\ref{t
large}.\QED



\bigskip

\noindent {\bf Alternative Proof of Theorem~\ref{coredim2}.} Again
assuming $\mu=1$ and $\delta=1$ we wish to prove that
$\core(I)=\m^{2n-1}$. But $\m^n$ is integral over $I$ and
$\core(\m^n)=\m^{2n-1}$ by Corollary 5.3 for instance. Hence
$\core(I) \supset \core(\m^n) =\m^{2n-1}$. Thus we only need to
establish the inclusion ${\rm core}(I) \subset {\mathfrak
m}^{2n-1}$.  Since $\core(I)$ is a monomial ideal it suffices to
prove that ${\mathfrak m}^{2n-1}$ is the maximal monomial ideal
contained in some reduction $J$ of $I$, i.e. ${\mathfrak
m}^{2n-1}=\mono(J)$. We take $J=(y^n-x^n,f)$ for $f = b_0y^n -
b_1x^{n-k_1}y^{k_1} - \cdots - b_sx^{n-k_s}y^{k_s}$ with $(b_0,
\ldots , b_s) \in k^{s+1}$ general. Notice
$\underline{\beta}=x^{2n},y^{2n}$ is a regular sequence of monomials
contained in $J$ and $(\bbb) \colon \m^{2n} =\m^{2n-1}$. Thus
according to Lemma~\ref{linkage} the equality $\mono(J)=\m^{2n-1}$
follows once we have shown that $\Mono((\bbb):J)= \m^{2n}$. To
compute  $(\bbb):J = (x^{2n},y^{2n}) \colon (y^n-x^n, f) $ we write
$x^{2n}= h (y^n-x^n) +gf$ where $h,g$ are forms of degree $n$ and
${\rm deg}_y  \, g \le n-1$. We have
\begin{eqnarray*}
x^{2n} & =& h (y^n-x^n) +gf \\
y^{2n} &=& (h + y^n + x^n) (y^n-x^n) + gf.
\end{eqnarray*}
Hence $(x^{2n},y^{2n}) \colon (y^n-x^n, f)= (x^{2n}, y^{2n},
\Delta), $ where
\[\Delta =
\left|%
\begin{array}{cc}
  h & g \\
  h+y^n +x^n & g \\
\end{array}%
\right|
= - (y^n+x^n) g.\] To prove that $\Mono(x^{2n}, y^{2n}, \Delta)=
\m^{2n}$ it suffices to show that the monomial support of
$\Delta=-(y^n+x^n) g$ is the set of all monomials of degree $2n$
except for $y^{2n}$. To this end we establish that the monomial
support of $g$ is the set of all monomials of degree $n$ except for
$y^n$. After dehomogenizing the latter claim follows from a general
fact about polynomials in $k[y]$:

\ms

\begin{Lemma}\label{k[y]} Let $k[y]$ be a polynomial ring
over an infinite field $k$, and $f = b_0y^n - b_1y^{k_1} - \ldots -
b_sy^{k_s} \in k[y],$ where   $0 < k_1 < \ldots < k_s < n$ are
integers with $\gcd(k_1, \ldots, k_s, n) = 1$ and $(b_0, \ldots ,
b_s) \in k^{s+1}$ is general. If $1 = h(y^n - 1) + gf$ with $h \in
k[y]$ and $g = c_0 + c_1y + \ldots + c_{n-1}y^{n-1} \in k[y]$, then
$c_i \ne 0 $ for every $i$.
\end{Lemma}

\medskip
\noindent To prove Lemma~\ref{k[y]} we are led to study Hankel
matrices with strings of zeros and variables. We need to determine
under which conditions on the distance between the strings of
variables the ideal generated by the maximal minors of the matrix
has generic grade. We solve this problem, which is interesting in
its own right, by using techniques from Gr\"obner basis theory. On
the other hand, Lemma 6.7 is actually equivalent to Theorem 6.4.
Therefore the first proof of Theorem 6.4 also provides a less
involved proof of Lemma 6.7.

\bigskip

\section{Monomials of the Same Degree: Dimension three}\label{secdim3}

\ms

In this section we study the core of ideals generated by monomials
of the same degree in three variables. However, our results are
less complete than in the two dimensional case.

\medskip

\begin{Notation and Discussion}\label{abc}
 {\rm Let $R = k[x,y,z]$ be a polynomial ring over an infinite
 field $k$ and consider the $R$-ideal
 $I =(x^n,y^n, z^n, \{x^{n - k_i }y^{k_i}\},
 \{x^{n -\ell_i}z^{\ell _i}\}, \{y^{n-m_i}z^{m_i}\} ) \not= R$.  Write
\begin{eqnarray*}
a &=& \gcd(n,k_i\mbox{'s},  \ell_i\mbox{'s} )  \\
 b &=& {\rm gcd}(n, k_i\mbox{'s}, m_i\mbox{'s} ) \\
 c &=& {\rm gcd}(n, \ell_i\mbox{'s},  m_i\mbox{'s} ) \\
 S &=& k[x^a, y^b, z^c].
\end{eqnarray*}
Notice that $\gcd(a,b) = \gcd(a,c) = \gcd(b,c) = \gcd(a,b,c)$. For
the purpose of computing the core of $I$ we may assume that $\delta=
\gcd(a,b,c) = 1$, since we may first compute the core of the
corresponding ideal in the polynomial ring $k[x^\delta, y^\delta,
z^\delta]$ and then use the fact that the core is preserved under
flat base change according to Proposition~\ref{local global} and
~\cite[4.8]{CPU1}. Thus throughout this section we will assume that
$gcd(a,b,c) = 1$, and hence that, $a$, $b$, $c$ are pairwise
relatively prime. Furthermore by relabeling the variables we can
assume that $ a \le b \le c$.

 Let $J$ be the $R$-ideal generated
by $x^n$, $y^n$, $z^n$, let $K$ be the $R$-ideal generated by the
monomials in $S$ of degree $n$, and $L$ the $R$-ideal generated by
$S_{\geq n}$. Clearly $J \subset I \subset K \subset L$. }
\end{Notation and Discussion}

\ms

We will show that the core of $I$ is always equal to the core of
$K$; in particular $K$ is contained in the first coefficient ideal
of $I$ according to Corollary~\ref{checkI3}(a). If $a=1$, we will
actually show that $\core(I)=\core(K)=\core(L)$ and that $L$ is the
first coefficient ideal of $I$. We first need some technical lemmas.
For their proofs set $k = \gcd(n, k_i\mbox{'s})$, $\ell = \gcd(n,
\ell_i\mbox{'s})$, and $m = \gcd(n, m_i\mbox{'s})$.

\ms

\begin{Lemma}\label{K is in}
With assumptions as in  \ref{abc} one has $K^{3t} \subset S_{3nt}
R \subset I^{3t} + J^{\langle t+1 \rangle}$ for $t \gg 0$.
\end{Lemma}
\demo   It suffices to show that for a  monomial
$x^{au}y^{bv}z^{cw}$ of $S_{3nt}$ that is not in $J^{\langle t+1
\rangle}$, we have $x^{au}y^{bv}z^{cw} \in I^{3t}$. Thus $au + bv +
cw = 3nt$ and  $au, bv, cw < n(t+1)$.  Since the sum of any two of
$au$, $bv$, $cw$ is strictly less than $2n(t+1) $ we have $au , bv,
cw > (t-2)n$.  In particular, when $t \gg 0$ each summand $au, bv,
cw \gg 0$.
Applying Lemma~\ref{numeric} to the integers $n$, $\ell_i$, $m_i$ we
can write
\begin{equation}\label{eqn cw}
cw = \alpha n + \sum \beta_i \ell_i + \sum \gamma_i m_i,
\end{equation}
where $\sum \beta_i + \sum \gamma_i < n/c$ and $\alpha, \beta_i,
\gamma_i \ge 0$. In particular
 \begin{equation}\label{eqn cw 1}
\alpha n = cw -( \sum \beta_i \ell_i + \sum \gamma_i m_i) > (t-2 -
n/c)n.
\end{equation}

Next we wish to apply Lemma~\ref{numeric} to the integers $n$,
$k_i$, $\ell m$.  Since $\sum \gamma_i(n - m_i) < n^2/c$ we have $bv
- \sum \gamma_i(n - m_i) \gg 0$. We first observe that $\gcd(n,
k_i\mbox{'s}, \ell m) = \gcd(k, \ell m)= ab$.  This follows since $a
= \gcd(k,\ell)$, $b = \gcd(k,m)$, and $\gcd(a,b) =1$. Now we want to
prove that  $bv - \sum \gamma_i(n - m_i)$ is divisible by $ab$.
Clearly $b$ divides $bv - \sum \gamma_i(n - m_i)$. Since $au = 3nt -
bv -\alpha n - \sum \beta_i \ell_i - \sum \gamma_i m_i$ by (\ref{eqn
cw}), we see that $a$ divides $bv + \sum \gamma_i m_i$ and hence
divides $bv - \sum \gamma_i(n - m_i)$. As $\gcd(a,b) =1$, $bv - \sum
\gamma_i(n - m_i)$ is a multiple of $ab$. Hence according to
Lemma~\ref{numeric} we can write
   \begin{equation}\label{eqn bv 0}
bv - \sum \gamma_i(n - m_i)  = \mu n + \sum \nu_i k_i + \eta \ell m,
    \end{equation}
        where  $\sum \nu_i + \eta < n/ab$
        and $\mu, \nu_i, \eta \ge 0$.
Therefore
 \begin{equation}\label{eqn bv 1}
   bv =   \mu n +  \sum \gamma_i(n - m_i)  + \sum \nu_i k_i + \eta \ell m.
   \end{equation}

   Now we apply Lemma~\ref{numeric} to the integers $n$, $n- m_i$.  By (\ref{eqn bv
   0})
   we have  $\mu n + \eta \ell m \gg 0$ as $\sum \nu_i k_i  < n^2/ab$. Hence we may write
   \begin{equation*}
\mu n + \eta \ell m = \rho n + \sum \gamma^{\, \prime}_i(n - m_i),
    \end{equation*}
    where $ \sum \gamma^{\, \prime}_i< n/m$ and $\rho,  \gamma^{\, \prime}_i   \ge 0$.
    Substituting the last equality into (\ref{eqn bv 1}) we obtain
    \begin{equation}\label{eqn bv 2}
  bv =   \rho n + \sum \gamma^{\, \prime}_i(n - m_i) + \sum \gamma_i(n - m_i)  + \sum \nu_i k_i .
    \end{equation}

Next consider $au - \sum \beta_i (n - \ell_i) - \sum \nu_i(n -
k_i)$,   which is $ \gg 0$ when $t \gg 0$.  We wish to see that $au
- \sum \beta_i (n - \ell_i) - \sum \nu_i(n - k_i)$ is divisible by
$\ell$. Indeed
\begin{eqnarray*}
   au   - \sum \beta_i (n - \ell_i) -
   \sum \nu_i(n - k_i)&\equiv &au +
    \sum \nu_i  k_i  \ \  {\rm mod} \ \ell \\
     & \equiv& au + cw - \sum  \gamma_i m_i +   \sum \nu_i  k_i  \ \  {\rm mod} \ \ell   \hspace{.45in} {\rm by} \ (\ref{eqn cw})\\
      &\equiv &  au + cw + bv    \ \  {\rm mod} \ \ell \hspace{1.28in} {\rm by} \ (\ref{eqn bv 1}) \\
      & \equiv& 3nt  \ \  {\rm mod} \ \ell \\
     & \equiv & 0  \ \  {\rm mod} \ \ell.
\end{eqnarray*}
Therefore  $au - \sum \beta_i (n - \ell_i) - \sum \nu_i(n - k_i)$ is
a multiple of  $\ell$.  Thus we may apply Lemma~\ref{numeric} to the
integers $n$, $n - \ell_i$ to write
\begin{equation*}
au - \sum \beta_i (n - \ell_i) - \sum \nu_i(n - k_i) = \zeta n +
\sum \beta^{\, \prime}_i(n - \ell_i),
\end{equation*}
where $ \sum \beta^{\, \prime}_i< n/\ell$ and $\zeta,  \beta^{\,
\prime}_i\ge 0$.   Hence
\begin{equation}\label{eqn au}
au =  \zeta n + \sum \beta_i (n - \ell_i) + \sum \nu_i(n - k_i)  +
\sum  \beta^{\, \prime}_i(n - \ell_i).
\end{equation}
Combining equations (\ref{eqn au}), (\ref{eqn bv 2}), and (
\ref{eqn cw}) we obtain
\begin{eqnarray*}
 (au, bv, cw)& =&  \zeta(n,0,0) +  \rho (0,n,0) + \alpha(0,0,n) +
      \sum ( \beta_i +  \beta^{\, \prime}_i)(n - \ell_i, 0 , \ell_i)   \\
    &+&   \sum \nu_i (n - k_i,  k_i , 0) +\sum (\gamma_i +   \gamma^{\, \prime}_i)( 0, n - m_i, m_i) \\
    & - & (0,0, \sum  \beta^{\, \prime}_i \ell_i  + \sum  \gamma^{\, \prime}_i  m_i) .
\end{eqnarray*}
      Taking the sum of the components on each side we see that
      $ \sum  \beta^{\, \prime}_i \ell_i  + \sum  \gamma^{\, \prime}_i  m_i = \lambda n$ for some $\lambda \ge 0$. Thus
\begin{eqnarray*}
 (au, bv, cw)& =&  \zeta(n,0,0) +  \rho (0,n,0) +( \alpha - \lambda)(0,0,n) +
      \sum ( \beta_i +  \beta^{\, \prime}_i)(n - \ell_i, 0 , \ell_i)   \\
    &+&   \sum \nu_i (n - k_i,  k_i , 0) +\sum (\gamma_i +   \gamma^{\, \prime}_i)( 0, n - m_i, m_i).
    \end{eqnarray*}

Since $ \sum \beta^{\, \prime}_i< n/\ell$ and $ \sum \gamma^{\,
\prime}_i< n/m$ we must have $\lambda n < (n/\ell + n/m)n$, and
consequently $\lambda < n/\ell + n/m$.  As $\alpha > t-2 - n/c$ by
(\ref{eqn cw 1}), we have $\alpha - \lambda \ge 0$ for $t \gg 0$.
Finally, since the sum of the components on the left hand side is
$3nt$ we deduce that the right hand side is the exponent vector of a
monomial in $I^{3t}$, as desired. \QED

\ms

\begin{Lemma}\label{gens}
With assumptions as in  \ref{abc} the $S$-ideal $S_{\ge j}$ is
generated by monomials of degrees at most $j+b-1$ for every integer
multiple $j$ of $c$.
\end{Lemma}
\demo Let  $x^{au}y^{bv}z^{cw}$ be a minimal monomial generator of
$S_{\ge j}$. Suppose that $au + bv + cw \ge j +b$. Since $a \leq b$
it follows that $u=v=0$ because the monomial $x^{au}y^{bv}z^{cw}$ is
a minimal generator of $S_{\ge j}$. Hence $cw \ge j+b
> j$ which implies $z^{cw}=z^{j}z^{c(w-j/c)}$, a contradiction.
\QED

\ms

\begin{Lemma}\label{K = L}
With assumptions as in \ref{abc} and $a=b=1$ the $S$-ideal $S_{\ge
j}$ is generated by monomials of degree $j$ for every integer
multiple $j$ of $c$; in particular $L=K$.
\end{Lemma}
\demo This follows immediately from Lemma~\ref{gens}. \QED

\ms

\begin{Lemma}\label{L is in}
With assumptions as in  \ref{abc} and $a=1$ one has $L^{3t}
\subset S_{\ge 3nt}R \subset I^{3t} + J^{\langle t+1 \rangle}$ for
$t \gg 0$.
\end{Lemma}
\demo
It suffices to show that every minimal monomial generator
$x^{u}y^{bv}z^{cw}$ of the $S$-ideal $S_{\ge 3nt}$ that is not in
$J^{\langle t+1 \rangle}$ is in $I^{3t}$. Lemma~\ref{gens} gives $u
+ bv + cw = 3nt + \epsilon$ with $0 \le \epsilon \le b-1$. Since
$x^{u}y^{bv}z^{cw} \not\in J^{\langle t+1 \rangle}$ we have  $ bv,
cw < n(t+1)$, hence $ bv + cw < 2n(t+1)$. As  $u + bv + cw \ge 3nt$
we obtain $u > (t-2)n$. In particular $u \ge \epsilon$ for $t \ge
3$. Now $x^{u}y^{bv}z^{cw}=x^{\epsilon}x^{u-\epsilon}y^{bv}z^{cw}$
with $x^{u-\epsilon}y^{bv}z^{cw} \in S_{3nt}R$, and the assertion
follows from Lemma~\ref{K is in}. \QED

\ms

From now on we will assume that the field $k$ is infinite.

\begin{Theorem}  With assumptions as in  \ref{abc} one has $\core(I) = \core(K)$.
In particular $K \subset \check{I}$, the first coefficient ideal
of $I$ .
\end{Theorem}
\demo Lemma~\ref{K is in} gives $K^{3t} +J^{\langle t+1 \rangle}=
I^{3t} + J^{\langle t+1 \rangle}$ for $t \gg 0$. Thus $\core(K)=
\core(I)$ by Theorem~\ref{core formula}. Corollary~\ref{checkI3}(a)
then implies that $\check{K} = \check{I}$.\QED

\ms

We are now ready to give an explicit formula for the core of $I$.

\begin{Theorem}\label{coreL}
With assumptions as in \ref{abc} and $a=1$ one has
\[{\rm core}(I)= \core(K) = \core(L) =  (S_{\geq 3n-b-c})R.
\]
\end{Theorem}
\demo The $R$-ideal $J=(x^n,y^n,z^n)$ is a reduction of $L$
according to Lemma~\ref{Sn}(a) and the $S$-ideal $S_{\ge n}$ is
normal by \cite[3.13]{Vi} and \cite[3.5]{RRV}. Now we obtain for $t
\gg 0$,
\begin{eqnarray*}
  J^{\langle t+1 \rangle} :_R L^{3t} &=& J^{t+1} :_R L^{t}  \hspace{.63in} \mbox{by Lemma~\ref{colon lemma}}\\
   &\subset& \core(L) \hspace{.74in} \mbox{by Proposition~\ref{local global} and \cite[4.8]{PU}}\\
   &\subset& \core(K) \hspace{.73in} \mbox{since $K$ is a reduction of $L$} \\
   &\subset& \core(I) \hspace{.78in} \mbox{since $I$ is a reduction of $K$}  \\
   &=& J^{\langle t+1 \rangle} :_R I^{3t}   \hspace{.55in} \mbox{by Theorem~\ref{core formula}}\\
   &=&  J^{\langle t+1 \rangle} :_R L^{3t} \hspace{.53in} \mbox{by Lemma~\ref{L is in}}\\
   &=&  (S_{\geq 3n-b-c})R \hspace{.48in} \mbox{by
   Lemma~\ref{Sn}(b)}.
\end{eqnarray*}
\QED

\ms

The next example shows that Theorem~\ref{coreL} does not hold when
$a=2$.

\begin{Example}
{\rm Let $R=k[x,y,z]$ be a polynomial ring over a field $k$ with
$\car \, k =0$ and consider the ideal $I=(x^{30}, y^{30}, z^{30},
x^{6}y^{24}, x^{10}z^{20}, y^{15}z^{15})$. In this case $a=2$,
$b=3$, $c=5$ and $S=k[x^2,y^3,z^5]$. One has $L = K + (x^{26}z^5,
x^{20}y^6z^5, x^{16}z^{15}, x^{14}y^{12}z^5, x^{10}y^{6}z^{15},
x^8y^{18}z^5, x^{4}y^{12}z^{15}, x^2y^{24}z^5) + (y^{27}z^5,
y^{12}z^{20})$. It turns out that $\core(L)=S_{\ge 81}R \subsetneq
\core(I)=\core(K)$.
}
\end{Example}

\ms

\begin{Theorem}\label{S2}
With assumptions as in \ref{abc} and $a=1$ one has
\begin{itemize}
\item[(a)] $\check{I}=L$. \item[(b)] $ R[\check{I}t]=R[Lt]$ is the
$S_2$-ification  of  $R[It]$.
\item[(c)] $ R[\check{I}t]=R[Lt]$ is a Cohen-Macaulay ring.
\end{itemize}
\end{Theorem}
\demo The ideal $L$ is integral over $I$ by Lemma~\ref{Sn}(a).
Furthermore $J^{t+1} \colon L^t=J^{t+1} \colon I^t$ for $t \gg 0$
according to Lemmas~\ref{L is in} and \ref{colon lemma}. Now
Corollary~\ref{checkI4} implies that $\check{L}=\check{I}$. Thus the
theorem follows once we have shown that $R[Lt]$ is Cohen-Macaulay.
The Rees algebra $S[S_{\ge n}t]$ is normal by
 \cite[3.13]{Vi} and \cite[3.5]{RRV}, and hence Cohen-Macaulay
 according to \cite[Theorem 1]{H}. But
$R[Lt]$ is a finite free module over $S[S_{\ge n}t]$ and thus a
Cohen-Macaulay ring as well. \QED

\ms

The next two corollaries show that for $a=b=1$ our formula for the
core becomes more explicit, akin to the case of two variables.

\begin{Corollary}\label{coreab=1}
In addition to the assumptions of \ref{abc} suppose that $a=b=1$
and write $q=\frac{3n}{c} -1$. One has
 \begin{itemize}
 \item[(a)]$\cI = K=L=((x,y)^{c}, z^{c})^{n/c}$.
 \item[(b)]$\core(I) = (z^{q c}) + \sum_{i=0}^{q-1} z^{i c}(x,y)^{(q-i)c
 -1}.$
 \end{itemize}
\end{Corollary}
\demo The first two equalities in part (a) follow from Lemma~\ref{K
= L} and Theorem~\ref{S2}(a), whereas the last equation is immediate
from the definition of $K$. To prove part (b) one uses
Theorem~\ref{coreL}.\QED

\ms

\begin{Corollary}
With assumptions as in \ref{abc} and  $a=b=c=1$ one has
 \begin{itemize}
 \item[(a)]$\cI =K=L=\m^n$.
 \item[(b)]$\core(I) =
 \m^{3n-2} = {\rm adj}(I^3).$
 \end{itemize}
\end{Corollary}
\demo In light of Corollary~\ref{coreab=1} it suffices to prove that
$\core(I)=\adj(I^3)$ in part (b). Indeed, part (a) and Lemma 7.2
show that the assumptions of Theorem 4.11 are satisfied. Now apply
that theorem. \QED

\end{document}